\documentclass[12pt, leqno]{amsart}
\usepackage{graphicx}
\usepackage{amsfonts,delarray,amssymb,amsmath, amsthm,a4,a4wide}
\usepackage{latexsym}
\usepackage{epsfig}

\newtheorem{thm}{Theorem}[section]
\newtheorem{cor}[thm]{Corollary}
\newtheorem{lem}[thm]{Lemma}
\newtheorem{prop}[thm]{Proposition}
\theoremstyle{definition}
\newtheorem{defn}[thm]{Definition}
\theoremstyle{remark}
\newtheorem{rem}[thm]{Remark}
\numberwithin{equation}{section}

\newcommand{\abs}[1]{\left\vert#1\right\vert}

\newcommand{\R}{\mathbb R}

\newcommand{\eps}{\epsilon}

\newcommand{\p}{\partial}

\newcommand{\comment}[1]{}
\def\h{\hspace*{.24in}} 
\newcommand*{\avint}{\mathop{\ooalign{$\int$\cr$-$}}}

\begin{document}

\title[Some minimization problems]{Some minimization problems in the class of convex functions with prescribed determinant}
\author{N. Q. Le}
\address{Department of Mathematics, Columbia University, New York, NY 10027}
\email{\tt  namle@math.columbia.edu}
\author{O. Savin}
\address{Department of Mathematics, Columbia University, New York, NY 10027}
\email{\tt  savin@math.columbia.edu}



\begin{abstract}
We consider minimizers of linear functionals of the type
$$L(u)=\int_{\p \Omega} u \, d \sigma - \int_{\Omega} u \, dx$$
in the class of convex functions $u$ with prescribed determinant $\det D^2 u =f$. 

We obtain compactness properties for such minimizers and discuss their regularity in two dimensions.
\end{abstract}
\maketitle

\section{Introduction}

In this paper, we consider minimizers of certain linear functionals in the class of convex functions with prescribed determinant. We 
are motivated by the study of convex minimizers $u$ for convex energies $E$ of the type
$$E(u) =\int_{\Omega} F(\det D^{2} u ) \, dx + L(u), 
\quad \quad \mbox{with $L$ a linear functional,}$$
which appear in the work of Donaldson \cite{D1}-\cite{D4} in the context of existence of K\"{a}hler metrics 
of constant scalar curvature for toric varieties. The minimizer $u$ solves a fourth order elliptic equation 
with two nonstandard boundary conditions involving the second and third order derivatives of $u$ (see \eqref{EL-intro} below). In this 
paper, we consider minimizers of $L$ (or $E$) in the case when the determinant $\det D^2u$ is prescribed. This allows us to understand better the type of boundary conditions that appear in such problems and to obtain estimates also for unconstrained minimizers of $E$.

The simplest minimization problem with prescribed determinant which is interesting in its own right is the following
$$\mbox{minimize} \, \,  \int_{\p \Omega} u \, d \sigma, \quad \mbox{ with $u \in \mathcal A_0,$}$$
where $\Omega$ is a bounded convex set, $d\sigma$ is the surface measure of $\p\Omega$, and $\mathcal A_0$ is the class 
of nonnegative solutions to the Monge-Amp\`ere equation $\det D^{2} u =1$:
$$\mathcal A_0:=\{ u: \bar \Omega \to [0, \infty)| \, \mbox{$u$ convex}, \quad \det D^2 u = 1\}.$$
{\it Question: Is the minimizer $u$ smooth up to the boundary $\p\Omega$ if $\Omega$ is a smooth, say uniformly convex, domain?}

In the present paper, we answer this question affirmatively in dimensions $n=2$. First, we remark 
that the minimizer must vanish at $x_{0}$, the center of mass of $\p\Omega$:
\begin{equation*}
 x_{0} =\avint _{\p\Omega} x \, d\sigma.
\end{equation*}
This follows easily since $$u(x) - u(x_0)-\nabla u (x_0) (x-x_0) \in \mathcal A_0$$
and
$$
\int_{\p\Omega} [u(x) - u(x_0)-\nabla u (x_0) (x-x_0)] d\sigma =\int_{\p \Omega} [u-u(x_0)] \, d \sigma 
\le \int_{\p \Omega} u \, d \sigma,
$$
with strict inequality if $u(x_0)>0$. Thus we can reformulate the problem above as minimizing $$\int_{\p \Omega} u \, d \sigma - \mathcal H^{n-1}(\p\Omega)\, \,  u(x_0)$$
in the set of all solutions to the Monge-Amp\`ere equation $\det D^2 u=1$ which are not necessarily nonnegative. This formulation is more convenient since now we can perturb functions in all directions.

 More generally, we consider linear functionals of the type
 \begin{equation*}
 L(u)= \int_{\p\Omega} u \, d\sigma -\int_{\Omega} u \, dA,
\label{L-fnc}
\end{equation*}
with $d\sigma$, $dA$ nonnegative Radon measures supported on $\p \Omega$ and $\Omega$ respectively. In this paper, we 
study the existence, uniqueness and regularity properties for minimizers of $L$. i.e.,
\begin{equation*}\label{Problem}
\mbox{(P)} \quad \quad \quad \mbox{minimize $L(u)$ for all $u \in \mathcal A$}
\end{equation*}
in the class $\mathcal A$ of subsolutions (solutions) to a Monge-Amp\`ere equation $\det D^{2} u \geq f$:
$$\mathcal A:=\{u: \overline \Omega \to \R| \, \mbox{$u$ convex}, \quad \det D^2 u \ge f \}.$$
Notice that we are minimizing a linear functional $L$ over a convex set $\mathcal A$ in the cone of convex functions.

Clearly, the minimizer of the problem (P) satisfies $\det D^2 u=f$ in $\Omega$. Otherwise we can find $v \in \mathcal A$ such that $v=u$ in a neighborhood of $\p \Omega$, and $v \ge u$ in $\Omega$ with strict inequality in some open subset, thus $L(v)<L(u)$.

We assume throughout that the following 5 conditions are satisfied:

1) $\Omega$ is a bounded, uniformly convex, $C^{1,1}$ domain.

2) $f$ is bounded away from $0$ and $\infty$.

3) $$ d\sigma = \sigma(x) \, d\mathcal{H}^{n-1}\lfloor \p\Omega,$$
with the density $\sigma(x)$ bounded away from $0$ and $\infty$.

4) $$dA=A(x) \, dx \quad \mbox{in a small neighborhhod of $\p \Omega$}$$
with the density $A(x)$ bounded from above.

5) $$L(u)>0 ~\text{for all} ~u ~\text{convex but not linear}.$$

\

The last condition is known as {\it the stability of $L$} (see \cite{D1}) and in 2D, is equivalent to saying that, for all 
linear functions $l$, we have $$L(l)=0 \quad \quad \mbox{and} \quad L(l^{+})>0 \quad \mbox{if} \quad l^{+}\not\equiv 0 \quad \mbox{in $\Omega$,}$$ where $l^{+}= \max (l, 0)$ (see Proposition \ref{stab_2d}).

Notice that the stability of $L$ implies that
$L(l) = 0$ for any linear function $l$,
hence $d \sigma$ and $d A$ must have the same mass and the same center of mass.

 A minimizer $u$ of the functional $L$ is determined up to linear functions since both $L$ and $\mathcal A$ are invariant 
under addition with linear functions. We ``normalize" $u$ by subtracting its the tangent plane at, say the center of 
mass of $\Omega$. In Section \ref{Sec2}, we shall prove in Proposition \ref{min-prop} that there exists a unique normalized minimizer to the problem (P). 
 
 We also prove a compactness theorem for minimizers.
\begin{thm}[Compactness]
Let $u_{k}$ be the normalized minimizers of the functionals $L_{k}$ with data $(f_{k}, d\sigma_{k}, dA_{k},\Omega)$
that has uniform bounds in $k$. Precisely, the inequalities \eqref{sbd} and \eqref{stab1} below are satisfied uniformly in $k$ and $\rho\leq f_{k}
\leq \rho^{-1}$. If
\begin{equation*}
 f_{k}\rightharpoonup f, ~d\sigma_{k}\rightharpoonup d\sigma,~ dA_{k}\rightharpoonup dA,
\end{equation*}
then $u_{k}\rightarrow u$ uniformly on compact sets of $\Omega$ where $u$ is the normalized minimizer of the functional $L$ with data $(f, d\sigma, dA,\Omega).$
\label{cpt}
\end{thm}

If $u$ is a minimizer, then the Euler-Lagrange equation reads (see Proposition \ref{EL})
$$ \mbox{if $\varphi:\Omega \to \R$ solves $U^{ij}\varphi_{ij}=0$  then   $L(\varphi)=0,$}$$
where $U^{ij}$ are the entries of the cofactor matrix $U$ of the Hessian $D^2 u$.
Since the linearized Monge-Amp\`ere equation is also an equation in divergence form, we can always express 
the $\Omega$-integral of a function $\varphi$ in terms of a boundary integral. For this, we consider the solution $v$ to the Dirichlet problem
$$ U^{ij} v_{ij}= -dA \quad \text{in}~ \Omega, \quad \quad
v =0\quad \text{on}~\p \Omega.\\\ $$
Integrating by parts twice and using $\p_i(U^{ij})=\p_j (U^{ij})=0$, we can compute
\begin{align}\label{intparts}
\nonumber \int_{\Omega} \varphi \, dA&=-\int_{\Omega}  \varphi \, U^{ij}v_{ij} \\
\nonumber &=\int_{\Omega} \varphi_i \, U^{ij} v_j \, -\int_{\p \Omega} \varphi U^{ij}v_j \nu_i  \\
&=-\int_{\Omega} (U^{ij} \varphi_{ij}) v \,  + \int_{\p \Omega} \varphi_i  U^{ij}v \nu_j -\int_{\p \Omega} \varphi U^{ij}v_j \nu_i \\
\nonumber &=-\int_{\p \Omega} \varphi \, \, U^{ij}v_i \nu_j .
\end{align}
From the Euler-Lagrange equation, we obtain $$ U^{ij}v_i \nu_j=- \sigma \quad \mbox{on $\p \Omega$}.$$ 
Since $v=0$ on $\p \Omega$, we have $v_i=v_\nu \nu_i$, and hence
$$U^{ij}v_i \nu_j=U^{ij}\nu_i \nu_j v _\nu =U^{\nu \nu} v_\nu=\left(\det D^2_{x'} u\right) v_{\nu}$$
with $x' \perp \nu$ denoting the tangential directions along $\p \Omega$.
 In conclusion, if $u$ is a smooth minimizer then there exists a function $v$ such that $(u,v)$ solves the system
\begin{equation}
 \left\{
 \begin{alignedat}{2}
   \det D^{2} u ~&=f \h~&&\text{in} ~\Omega, \\\
 U^{ij} v_{ij}&= -dA \h~&&\text{in}~ \Omega,\\\
v &=0\h~&&\text{on}~\p \Omega,\\\
U^{\nu \nu} v_{\nu} &= -\sigma~&&\text{on}~\p \Omega.
 \end{alignedat}
 \right.
\label{EL-eq-intro}
\end{equation}

 This system is interesting since the function $v$ above satisfies two boundary conditions, Dirichlet and Neumann, while $u$ has no boundary conditions. Heuristically, the boundary values for $u$ can be recovered from the term $U^{\nu \nu}=\det D^2_{x'}u$ which appears in the Neumann boundary condition
 for $v$.

Our main regularity results for the minimizers $u$ are in two dimensions.

\begin{thm}
 Assume that $n=2$, and the conditions 1)-5) hold. If $\sigma\in C^{\alpha}(\p\Omega)$,   $f \in C^{\alpha}(\overline{\Omega}),$ and
$\p \Omega \in C^{2,\alpha},$ then the minimizer $u\in C^{2,\alpha}(\overline{\Omega})$ and the system \eqref{EL-eq-intro} holds in the classical sense.
\label{C2-thm}
\end{thm}

We obtain Theorem \ref{C2-thm} by showing that $u$ separates quadratically on $\p \Omega$ from its tangent planes and then 
we apply the boundary H\"older gradient estimates for $v$ which were obtained in \cite{LS}.

As a consequence of Theorem \ref{C2-thm}, we obtain higher regularity if the data $(f,d\sigma, dA,\Omega)$ is more regular.
\begin{thm}
 Assume that $n=2$ and the conditions 1)-5) hold. If $\sigma\in C^{\infty}(\p\Omega)$, $f \in C^{\infty}(\overline{\Omega})$, $A\in C^{\infty}(\overline{\Omega})$, $\p \Omega \in C^{\infty},$ then $u\in C^{\infty}(\overline{\Omega})$.
\label{smooth-thm}
\end{thm}

In Section \ref{Sec-sing}, we provide an example of Pogorelov type for a minimizer in dimensions $n \ge 3$ that shows that Theorem \ref{smooth-thm} does not hold in this generality in higher dimensions. 

We explain briefly how Theorem \ref{smooth-thm} follows from Theorem \ref{C2-thm}. If $u \in C^{2,\alpha}(\overline \Omega)$, then $U^{ij} \in C^{\alpha}(\overline \Omega)$ and Schauder estimates give $v \in C^{2, \alpha}(\overline \Omega)$, thus $v_\nu \in C^{1,\alpha}(\p \Omega)$. From the last equation in \eqref{EL-eq-intro} we obtain $U^{\nu \nu}=\det D^2_{x'}u \in C^{1,\alpha}(\p \Omega)$. This implies $u \in C^{3,\alpha}(\p \Omega)$ and from the first equation in  \eqref{EL-eq-intro} we find $u \in C^{3,\alpha}(\overline \Omega)$. We can repeat the same argument and obtain that $u\in C^{k, \alpha}$ for any $k \ge 2$.

  As we mentioned above, our constraint minimization problem is motivated by the minimization of the Mabuchi energy functional from complex geometry in the case of toric varieties
\begin{equation*}
M(u) =\int_{\Omega} -\log \det D^{2} u + \int_{\p \Omega} u d\sigma -\int_{\Omega} u dA.
 \label{Mfunc}
\end{equation*}
In this case, $d\sigma$ and $dA$ are canonical measures on $\p\Omega$ and $\Omega$. Minimizers of $M$ satisfy the following fourth order equation,
called Abreu's equation \cite{Abreu}
\begin{equation*}
u^{ij}_{ij}:= \sum_{i, j=1}^{n}\frac{\p ^2 u^{ij}}{\p x_{i}\p x_{j}}= -A,
\label{Abreu}
\end{equation*}
where $u^{ij}$ are the entries of the inverse matrix of $D^2u$.
This equation and the functional $M$ have been studied extensively by Donaldson in a series of papers \cite{D1}-\cite{D4} (see also \cite{ZZ}). In these papers, the domain $\Omega$
was taken to be a polytope $P\subset \R^{n}$ and $A$ was taken to be a positive constant. The existence of smooth solutions with suitable boundary
conditions has important implications in complex geometry. It says that we can find K\"{a}hler metrics of constant scalar curvature
for toric varieties.

More generally, one can consider minimizers of the following convex functional
\begin{equation}\label{E}
 E(u) =\int_{\Omega} F(\det D^{2} u ) + \int_{\p \Omega} ud\sigma -\int_{\Omega} u dA
\end{equation}
 where $F(t^n)$ is a convex and decreasing function of $t \ge 0$. The Mabuchi energy functional 
corresponds to $F (t) =-\log t$ whereas in our minimization problem (P) (with $f\equiv 1$)
\begin{equation*}
 F (t) = \left\{
 \begin{alignedat}{1}
 \infty   \h~&\text{if}~ t<1, \\\
 0 \h~&\text{if}~t\geq 1.
 \end{alignedat}
 \right.
\end{equation*}
Minimizers of $E$ satisfy a system similar to \eqref{EL-eq-intro}:
\begin{equation}\label{EL-intro}
 \left\{
 \begin{alignedat}{2}
   -F'(\det D^{2} u) ~&=v \h~&&\text{in} ~\Omega, \\\
 U^{ij} v_{ij}&= -dA \h~&&\text{in}~ \Omega,\\\
v &=0\h~&&\text{on}~\p \Omega,\\\
U^{\nu \nu} v_{\nu} &= -\sigma~&&\text{on}~\p \Omega.
 \end{alignedat}
 \right.
\end{equation}

A similar system but with different boundary conditions was investigated by Trudinger and Wang in \cite{TW2}.
If the function $F$ is strictly decreasing then we see from the first and third 
equations above that $\det D^2 u=\infty$ on $\p \Omega$, and therefore we cannot expect 
minimizers to be smooth up to the boundary (as is the case with the Mabuchi functional $M(u)$).

If $F$ is constant for large values of $t$ (as in the case we considered) then $\det D^2 u$ becomes finite on the boundary and smoothness up to the boundary is expected. More precisely assume that
$$F \in C^{1,1}((0, \infty)), \quad G(t):= F(t^n) \quad \mbox{is convex in $t$,} \quad \mbox{and} \quad G'(0^+)=-\infty,$$
and there exists $t_0>0$ such that 
$$F(t)=0 \quad \mbox{on $[t_0, \infty)$}, \quad \quad F''(t)>0  \quad\mbox{on $(0,t_0]$.} $$

\begin{thm}\label{min_func}
Assume $n=2$, and the conditions 1)-5) and the above hypotheses on $F$ are satisfied. If  
$\sigma\in C^{\alpha} (\p\Omega), \, A \, \in C^\alpha(\overline \Omega),\,\p\Omega\in C^{2,\alpha}$ then the 
normalized minimizer $u$ of the functional $E$ defined in \eqref{E} satisfies $u\in C^{2,\alpha}(\overline \Omega)$ and the system \eqref{EL-intro} holds in the classical sense.  
\end{thm}

The paper is organized as follows. In Section 2, we discuss the notion of stability for the 
functional $L$ and prove existence, uniqueness and compactness of minimizers of the problem (P). In Section 3, we 
state a quantitative version of Theorem \ref{C2-thm}, Proposition \ref{singular}, and we also obtain the 
Euler-Lagrange equation. Proposition \ref{singular} is proved in sections 4 and 5, first under the 
assumption that the density $A$ is bounded from below and then in the general case. In Section \ref{Sec-sing}, we 
give an example of a singular minimizer in dimension $n \ge 3$. Finally, in Section 7, we prove Theorem \ref{min_func}. 

\section{Stability inequality and existence of minimizers}\label{Sec2}

Let $\Omega$ be a bounded convex set and define
 $$L(u)=\int_{\p \Omega} u \, d \sigma-\int_\Omega u \, d A$$
for all convex functions $u:\overline{\Omega} \to \R$ with $u \in L^1(\p \Omega, d \sigma)$.
We assume that
 \begin{equation}\label{sbd}
 \mbox{$\sigma \ge \rho$ on $\p \Omega$ and $A(x)\le \rho^{-1}$ in a neighborhood of $\p \Omega$,}
 \end{equation}
 for some small $\rho>0,$ and that $L$ is stable, i.e.,
\begin{equation}
L(u)>0 ~\text{for all} ~u ~\text{convex but not linear}.
\label{Lcon}
\end{equation}

Assume for simplicity that $0$ is the center of mass of $\Omega$. We notice that \eqref{Lcon} implies $L(l)=0$ for 
any $l$ linear since $l$ can be approximated by both convex and concave functions. We ``normalize" a convex function 
by subtracting its tangent plane at $0$, and this does not change the value of $L$. First, we prove some lower semicontinuity properties of $L$ with respect to normalized solutions.

 \begin{lem}[Lower semicontinuity]\label{l_semi}
  Assume that \eqref{sbd} holds and $(u_k)$ is a normalized sequence that satisfies
 \begin{equation}\label{int_bd}
 \int_{\p \Omega} u_k \, d\sigma \le C, \quad \quad u_k \to u \quad \mbox{uniformly on compact sets of $\Omega$,}
 \end{equation}
   for some function $u:\Omega \to \R$. Let $\bar u$ be the minimal convex extension of $u$ to $\overline \Omega$, i.e.,
 $$\bar u=u \quad \mbox{in $\Omega$}, \quad \bar u(x)=\lim_{t \to 1^-}u(tx) \quad \mbox{if $x \in \p \Omega$.}$$
 Then $$\int_\Omega u \, dA=\lim \int_\Omega u_k \, d A, \quad \quad \int_{\p \Omega} \bar u \, d\sigma \le \liminf 
\int_{\p \Omega} u_k \, d\sigma,$$ and thus
 $$ L(\bar u) \le \liminf L(u_k).$$
 \end{lem}

{\it Remark:} The function $\bar u$ has the property that its upper graph is the closure of the upper graph of $u$ in $\R^{n+1}$.

 \begin{proof}
 Since $u_k$ are normalized, they are increasing on each ray out of the origin. For each $\eta>0$ small, we consider
the set $\Omega_{\eta} :=\{x\in\Omega: dist(x,\p\Omega) < \eta\}$, and from \eqref{sbd} we obtain
\begin{equation*}
 \int_{\Omega_{\eta}} u_k \, dA\leq C \rho^{-1}\eta  \int_{\p\Omega} u_k \, d\sigma \le C \eta.
\end{equation*}
 Since this inequality holds for all small $\eta \to 0$, we easily obtain $$\int_\Omega u \, dA=\lim \int_\Omega u_k \, d A.$$

For each $z \in \p \Omega$, and $t<1$ we have $u_k(tz) \le u_k(z)$. We let $k \to \infty$ in the inequality
$$\int_{\p \Omega} u_k(tz) \, d \sigma \le \int_{\p \Omega} u_k(z) \, d \sigma$$ and obtain
$$\int_{\p \Omega} u(tz) \, d \sigma \le \liminf \int_{\p \Omega} u_k(z) \, d \sigma,$$
and then we let $t \to 1^-$,
$$\int_{\p \Omega} \bar u \, d \sigma \le \liminf \int_{\p \Omega} u_k \, d \sigma.$$
 \end{proof}

 \begin{rem}\label{weak_conv}From the proof we see that if we are given functionals $L_k$ with measures $\sigma_k$, $A_k$ that satisfy \eqref{sbd} uniformly in $k$ and $$\sigma_k \rightharpoonup \sigma , \quad A_k \rightharpoonup A,$$
and if \eqref{int_bd} holds for a sequence $u_k$, then the statement still holds, i.e.,
$$L(\bar{u}) \le \liminf L_k(u_k).$$
\end{rem}

 By compactness, one can obtain a quantitative version of \eqref{Lcon} known as {\it stablity inequality}. This was done 
by Donaldson, see Proposition 5.2.2 in \cite{D1}. For completeness, we sketch its proof here.

 \begin{prop}
Assume that (\ref{sbd}) and (\ref{Lcon}) hold. Then we can find $\mu>0$ such that
\begin{equation}
 L(u): =\int_{\p \Omega} ud\sigma -\int_{\Omega} u dA\geq \mu \int_{\p \Omega} u d\sigma
\label{stab1}
\end{equation}
for all convex functions $u$ normalized at $0$.
\label{stabi-prop}
\end{prop}
\begin{proof}
Assume the conclusion does not hold, so there is a sequence of normalized
convex functions $(u_{k})$ with
\begin{equation*}
 \int_{\p \Omega} u_{k} d\sigma =1, \quad \lim L(u_k)=0,
\end{equation*}
thus
\begin{equation*}
 \lim \int_{ \Omega} u_{n} dA = 1.
\end{equation*}
Using convexity, we may assume that $u_{k}$ converges
 uniformly on compact
subsets of $\Omega$ to a limiting function $u\ge 0$. Let $\bar{u}$ be the minimal convex
extension of $u$ to $\overline{\Omega}$. Then, from Lemma \ref{l_semi}, we obtain
$$L(\bar u)=0, \quad \int_\Omega \bar u \, dA=1,$$
thus $\bar u \ge 0$ is not linear, and we contradict \eqref{Lcon}.
\end{proof}

Donaldson showed that when $n=2$, the stability condition can be checked easily (see Proposition  5.3.1 in \cite{D1}).

\begin{prop}\label{stab_2d}
Assume $n=2$, \eqref{sbd} holds and
for all linear functions $l$ we have
\begin{equation}\label{st_2d}
L(l)=0 \quad \quad \mbox{and} \quad L(l^{+})>0 \quad \mbox{if} \quad l^{+}\not\equiv 0 \quad \mbox{in $\Omega$,}
\end{equation}
 where $l^{+}= \max (l, 0)$. Then $L$ is stable, i.e., condition \eqref{Lcon} is satisfied.
\end{prop}

\begin{proof}
For completeness, we sketch the proof. Assume by contradiction that $L(u)\le 0$ for some convex function $u$ 
which is not linear in $\Omega$. Let $u^*$ be the convex envelope generated by the boundary 
values of $\bar u$ - the minimal convex
extension of $u$ to $\overline{\Omega}$. Notice that $u^*=\bar u$ on $\p \Omega$. Since 
$L(u^*) \le L(\bar u) \le L(u)$ we find $L(u^*) \le 0$. Notice that $u^*$ is not linear since 
otherwise $0=L(u^*)< L(\bar u) \le 0$ (we used that $\bar u$ is not linear). After subtracting a linear function we 
may assume that $u^*$ is normalized and $u^*$ is not identically $0$.

We obtain a contradiction by showing that $u^*$ satisfies the stability inequality. By our hypotheses there exists $\mu>0$ small such that
$$L(l^+) \ge \mu \int_{\p \Omega} l^+ \, d \sigma,$$
for any $l^+$. Indeed, by \eqref{sbd} this inequality is valid if the ``crease" $\{l=0\}$ is 
near $\p \Omega$ and for all other $l$'s, it follows by compactness from \eqref{st_2d}. We approximate 
from below $u^*$ by $u^*_k$ which is defined as the maximum of the tangent planes of $u^*$ at some 
points $y_i \in \Omega$, $i=1,..,k$. Since $u^*$ is a convex envelope in 2D, $u^*_k$ is a discrete sum of $l^+$'s hence 
it satisfies the stability inequality. Now we let $k \to \infty$; since $u^*_k \le u^*$, using Lemma \ref{l_semi}, we obtain 
that $u^*$ also satisfies the stability inequality.
\end{proof}

\begin{prop}
 Assume that (\ref{sbd}) and (\ref{Lcon}) hold. Then there exists a unique (up to linear functions) minimizer $u$ of $L$ subject to the constraint $$u \in \mathcal A:=\{v: \overline \Omega \to \R| \, \mbox{$v$ convex}, \quad \det D^2 v \ge f \},$$
 where $\rho \le f \le \rho^{-1}$ for some $\rho >0$. Moreover, $\det D^2 u=f$.
\label{min-prop}
\end{prop}

\begin{proof}
Let $(u_k)$ be a sequence of normalized solutions such that $L(u_k) \to \inf_{\mathcal A} L$. By the stability inequality, we 
see that $\int_{\p \Omega} u_k \, d\sigma$ are uniformly bounded, and after passing to a subsequence, we may assume 
that $u_k$ converges uniformly on compact subsets of $\Omega$ to a function $u$.
Then $u \in \mathcal A$ and from the lower semicontinuity, we see that $L(u)=\inf_\mathcal A L$, i.e., $u$ is a 
minimizer. Notice that $\det D^2u = f$. Indeed, if a quadratic 
polynomial $P$ with $\det D^2 P>f$ touches $u$ strictly by below at some point $x_0 \in \Omega$, in a neighborhood of $x_0$, then we can replace $u$ in this neighborhood by $\max \{P+\eps,u\}\in \mathcal A$, and the energy decreases.

 Next we assume $w$ is another minimizer. We use the strict concavity of $M \mapsto \log (\det D^2 M)$ in the space of 
positive symmetric matrices $M$, and obtain that for a.e. $x$ where $u$, $w$ are twice differentiable
 $$\log \det D^2(\frac {u+w}{2})(x) \ge \frac 12 \log \det D^2 u (x) + \frac 12 \log \det D^2 w (x) \ge \log f(x).$$
  This implies $(u+w)/2 \in \mathcal A$ is also a minimizer and $D^2u=D^2w$ a.e in $\Omega$. Since $f$ is bounded above and below we know that $u,w \in W_{loc}^{2,1}$ (see \cite{DF}) in the open set $\Omega'$ where both $u,w$ are strictly 
convex. This gives that $u-w$ is linear on each connected component of $\Omega'$. If $n=2$, then $\Omega'=\Omega$ hence
 $u-w$ is linear. If $n\ge 3$, Labutin showed in \cite{L} that the closed set $\Omega \setminus \Omega'$ has 
Hausdorff dimension $n-2+2/n<n-1$, hence $\Omega'$ is connected, and we obtain the same conclusion 
that $u-w$ is linear in $\Omega$.\end{proof}


{\it Remark}: The arguments above show that the stability condition is also necessary 
for the existence of a minimizer. Indeed, if $u$ is a minimizer and $L(u_0)=0$ for some convex function $u_0$ that 
is not linear, then $u+u_0$ is also a minimizer and we contradict the uniqueness.

\

\begin{proof}[Proof of Theorem \ref{cpt}] We assume that the data $(f_{k}, d\sigma_{k}, dA_{k},\Omega)$ satisfies
 \eqref{sbd}, \eqref{stab1} uniformly in $k$ and $\rho \le f_k \le \rho^{-1}$. For each $k$, let $w_{k}$ be the 
convex solution to $\det D^{2} w_{k} = f_{k}$ in $\Omega$ with $w_{k} = 0$ on $\p\Omega$. Since $f_k$ are bounded 
from above we find $w_{k}\geq -C$,
and so by the minimality of $u_{k}$
\begin{equation*}
 L_{k} (u_{k}) \leq L_{k} (w_{k}) \le C.
\end{equation*}
It follows from the stability inequality that
\begin{equation*}
 \int_{\p\Omega} u_{k}\, d\sigma_{k} \leq C,
\end{equation*}
and we may assume, after passing to a subsequence, that $u_k \to u$ uniformly on compact sets.

We need to show that $u$ is a minimizer for $L$ with data $(f, d\sigma, dA,\Omega)$. For this it suffices to prove 
that for any continuous $v:\overline \Omega \to \R$ which solves $\det D^2 v =f$ in $\Omega$, we have
$L(u) \le L(v).$

Let $v_k$ be the solution to $\det D^{2} v_{k} = f_{k}$ with boundary data $v_{k} = v$ on $\p\Omega$. Using appropriate barriers it is standard to check that $f_k \rightharpoonup f$, $f_k \le \rho^{-1}$ implies $v_k \to v$ uniformly in $\overline \Omega$. Then, we let $k \to \infty$ in $L_k(u_k) \le L_k(v_k)$, use Remark \ref{weak_conv} and obtain
$$L(u) \le \liminf L_k(u_k) \le \lim L_k(v_k)=L(v),$$
which finishes the proof.
\end{proof}

\section{Preliminaries and the Euler-Lagrange equation}

We rewrite our main hypotheses in a quantitative way. We assume that for some small $\rho>0$ we have

H1) the curvatures of $\p \Omega$ are bounded from below by $\rho$ and from above by $\rho^{-1}$;

H2) $\rho \le f \le \rho^{-1};$

H3) $ d\sigma = \sigma(x) \, d\mathcal{H}^{n-1}\lfloor \p\Omega,$ with  $\rho\leq \sigma(x) \leq \rho^{-1};$

H4) $dA=A(x) \, dx$ in a small neighborhood $\Omega_\rho:=\{x \in \Omega| \quad dist \, (x, \p \Omega) < \rho \}$ of $\p \Omega$
with $A(x)\leq \rho^{-1}.$

H5) for any convex function $u$ normalized at the center of mass of $\Omega$, we have $$L(u):=\int_{\p \Omega} u \, d \sigma-\int_\Omega u \, dA \ge \rho \int_{\p \Omega} u \, d\sigma.$$

We denote by $c$, $C$ positive constants depending on $\rho$, and their values may change from line to line whenever
there is no possibility of confusion. We refer to such constants as {\it universal constants}.

Our main theorem, Theorem \ref{C2-thm}, follows from the next proposition which deals with less regular data.
\begin{prop}
 Assume that $n=2$ and the conditions H1-H5 hold. \\
(i) Then the minimizer $u$ obtained in Proposition \ref{min-prop} satisfies $u\in 
C^{1,\beta}(\overline \Omega) \cap C^{1,1}(\p\Omega)$ for some universal $\beta\in (0,1)$ and $u$ separates quadratically from its
tangent planes on $\p\Omega$, i.e.,
$$C^{-1}|x-y|^2 \le u(y)-u(x)-\nabla u(x)(y-x) \le C|x-y|^2 , \quad \quad \forall x,y \in \p \Omega,$$
for some $C>0$ universal.
\\
(ii) If in addition $\sigma\in C^{\alpha}(\p\Omega)$, then $u\mid_{\p\Omega}\in C^{2,\gamma}(\p\Omega),$  with $\gamma:=\min\{\alpha, \beta\}$ and
$$\|u\|_{ C^{2,\gamma}(\p\Omega)} \le C \|\sigma\|_{C^\gamma(\p \Omega)}.$$
\label{singular}
\end{prop}

It is interesting to remark that in part (ii), we obtain $u \in C^{2, \gamma}(\p \Omega)$ even though $f$ and $A$ are assumed to be only $L^{\infty}$.

\begin{proof}[Proposition \ref{singular} implies Theorem \ref{C2-thm}]

Theorem 7.3 in \cite{S2} states that a solution to the Monge-Amp\`ere equation which 
separates quadratically from its tangent planes on the boundary satisfies the classical $C^\alpha$-Schauder 
estimates. Thus, if the assumptions of Proposition \ref{singular} ii) are satisfied 
and $f \in C^\alpha(\overline \Omega)$ then $u\in C^{2,\gamma}(\overline{\Omega})$ with its $C^{2,\gamma}$ norm 
bounded by a constant $C$ depending on $\rho$, $\alpha$, $\| \sigma \|_{C^\alpha(\p \Omega)}$,
$\|\p\Omega\|_{C^{2,\alpha}}$, and $\|f\|_{C^\alpha(\overline \Omega)}$. This implies that the system \eqref{EL-eq-intro} holds in the classical sense. If $\alpha \le \beta$ then we are done. If $\alpha > \beta$ then we use $v\in C^{2, \beta}(\overline \Omega)$ in the last equation of the system and obtain $u \in C^{2, \alpha}(\p \Omega)$ which gives $u \in C^{2,\alpha}(\overline \Omega)$.
\end{proof}

We prove Proposition \ref{singular} in the next two sections. Part (ii) follows from part (i) and the boundary 
Harnack inequality for the linearized Monge-Amp\`ere equation which was obtained in \cite{LS} (see Theorem 2.4). This 
theorem states that if a solution to the Monge-Amp\`ere equation with bounded right hand side separates quadratically from its tangent planes on the boundary, then the classical boundary  estimate of Krylov holds for solutions of the associated linearized equation.

In order to simplify the ideas we prove the proposition in the case when the hypotheses H1, H2, H4 are replaced by

H1') $\Omega=B_1$;

H2') $f\in C^\infty(\overline \Omega)$, $\rho \le f \le \rho^{-1}$;

H4') $dA=A(x) \, dx$ with $\rho \le A(x) \le \rho^{-1}$ in $\Omega$ and $A \in C^\infty(\Omega)$.

We use H1' only for simplicity of notation. We will see from the proofs that the same arguments 
carry to the general case. We use H2' so that $D^2 u$ is continuous in $\Omega$ and the linearized Monge-Ampere equation 
is well defined. Our estimates do not depend on the smoothness of $f$, thus the general case 
follows by approximation from Theorem \ref{cpt}. Later in section \ref{5} we show that H4' can be replaced by H4, i.e the bounds for $A$ from below and above are not needed.

First, we establish a result on uniform modulus of convexity for minimizers of $L$ in 2D. 

\begin{prop}
Let $u$ be a minimizer of $L$ that satisfies the hypotheses above. 
Then, for any $\delta<1$, there exist $c(\delta)>0$ depending 
on $\rho$, $\delta$ such that $$ x \in B_{1-\delta} \Rightarrow \quad S_h(x) \subset \subset B_1 \quad \mbox {if} \quad h \le c(\delta).$$
\label{strict-prop}
\end{prop}
In the above proposition, we denoted by $S_{h}(x)$ the section of $u$ centered at $x$ at height $h$:
\begin{equation*}
 S_{h}(x)=\{y\in \bar{B_{1}}: u(y)< u(x) + \nabla u(x)(y-x) + h\}.
\end{equation*}

This result is well-known (see, e.g., Remark 3.2 in \cite{TW3}). For completeness, we include its proof here.
\begin{proof}
Without loss of generality assume $u$ is normalized in $B_1$, that is $u\geq 0$, $u(0) =0$.
From the stability inequality (\ref{stab1}), we obtain
\begin{equation*}
 \int_{\p B_1} u dx \leq C.
\end{equation*}
 This integral bound and the convexity of $u$ imply
\begin{equation*}
 |u|,\, \abs{Du}\leq C(\delta)~\text{in}~ B_{1-\delta/2},
\end{equation*}
for any $\delta<1.$ We show that our statement follows from these bounds.
Assume by contradiction that the conclusion is not true. Then, we can find a sequence of convex functions $u_{k}$ satisfying the bounds above such that
\begin{equation}
u_{k}(y_{k})\leq u_{k}(x_k) + \nabla u_{k}(x_k) (y_k-x_k) + h_{k}
 \label{uk}
\end{equation}
for sequences $x_{k}\in B_{1-\delta}$, $y_k \in \p B_{1-\delta/2}$ and $h_{k}\rightarrow 0.$
Because
$Du_{k}$
is uniformly bounded, after passing to a subsequence if necessary, we may assume
$$u_{k} \to u_* \quad \mbox{uniformly on $\overline B_{1-\delta/2} $}, \quad x_k \to x_*, \quad y_k \to y_*.$$
Moreover $u_*$ satisfies $\rho \le \det D^{2} u_* \le \rho^{-1},$ and
\begin{equation*}
 u_*(y_*) = u_*(x_*) + \nabla u_*(x_*) (y_*-x_*),
\end{equation*}
i.e., the graph of $u_*$ contains a straight-line in the interior.
However, any subsolution $v$ to $\det D^{2} v \ge \rho$ in $2D$ does
not have this property and we reached a contradiction.
\end{proof}

Since $f \in C^\alpha$ we obtain that $u \in C^{2,\alpha}(B_1)$ thus the linearized Monge-Amp\`ere equation is well defined in $ B_1$. Next lemma deals with general linear elliptic equations in $B_1$ which may become degenerate as we approach $\p B_1$.

\begin{lem}\label{linear}
Let $\mathcal L v:=a^{ij}(x) v_{ij}$
be a linear elliptic operator with continuous coefficients $a^{ij} \in C^\alpha(B_1)$ that satisfy the ellipticity condition $(a^{ij}(x))_{ij}>0$ in $B_1$. Given a continuous boundary data $\varphi$, there exists a unique solution $v \in C(\overline B_1) \cap C^2(\Omega)$ to the Dirichlet problem
$$\mathcal L v=0\quad \mbox{in $B_1$}, \quad v=\varphi \quad \mbox{on $\p B_1$.}$$
\end{lem}

\begin{proof}
For each small $\delta$, we consider the standard Dirichlet problem for uniformly elliptic 
equations $\mathcal L v_\delta=0$ in $B_{1-\delta}$, $v_\delta=\varphi$ on $\p B_{1-\delta}$. Since $v_\delta$ satisfies 
the comparison principle with linear functions, it follows that the modulus of continuity of $v_\delta$ at points on the 
boundary $\p B_{1-\delta}$ depends only on the modulus of continuity of $\varphi$. Thus, from maximum principle, we see that $v_\delta$ converges uniformly to a solution $v$ of the Dirichlet problem above.
The uniqueness of $v$ follows from the standard comparison principle.
\end{proof}

\begin{rem}
The modulus of continuity of $v$ at points on $\p B_1$ depends only on the modulus of continuity of $\varphi$.
\end{rem}
\begin{rem}\label{seq}
If $\mathcal L_m$ is a sequence of operators satisfying the hypotheses 
of Lemma \ref{linear} with $a^{ij}_m \to a^{ij}$ uniformly on compact subsets of $B_1$ and $\mathcal L_m v_m=0$ in $B_1$, $v_m=\varphi$ on $\p B_1$, then $v_m \to v$ uniformly in $\overline B_1$.

Indeed, since $v_m$ have a uniform modulus of continuity on $\p B_1$ and, for all large $m$, a uniform modulus of continuity in any ball $B_{1-\delta}$, we see that we can always extract a uniform convergent subsequence in $\overline B_1$. Now it is straightforward to check that the limiting function $v$ satisfies $\mathcal L v=0$ in the viscosity sense.
\end{rem}

Next, we establish an integral form of the Euler-Lagrange equations for the minimizers of L.
\begin{prop} \label{EL}
Assume that $u$ is the normalized minimizer of L in the class $\mathcal A$.
 If $\varphi\in C^{2}(\Omega)\cap C^{0}(\overline{\Omega})$ is a solution to the linearized Monge-Amp\`ere equation
$$U^{ij} \varphi_{ij} =0 \h~\text{in} ~\Omega,$$
 then
\begin{equation*}
 L(\varphi):= \int_{\p \Omega} \varphi d\sigma -\int_{\Omega} \varphi dA =0.
\end{equation*}
\label{lem-var}
\end{prop}

\begin{proof}
Consider the solution $u_{\eps} = u + \eps \varphi_{\eps}$ to
\begin{equation*}
 \left\{
 \begin{alignedat}{2}
   \det D^{2} u_{\eps} ~&=f \h~&&\text{in} ~B_1, \\\
u_{\eps} &= u + \eps \varphi\h~&&\text{on}~\p B_1.
 \end{alignedat}
  \right.
\end{equation*}
Since $\varphi_\eps$ satisfies comparison principle and comparison with planes, its existence follows as in Lemma \ref{linear} by solving the Dirichlet problems in $B_{1-\delta}$ and then letting $\delta \to 0$.

In $B_1$, $\varphi_{\eps}$ satisfies
\begin{equation*}
 0=\frac{1}{\eps}(\det D^{2} u_{\eps}-\det D^{2} u) = 
\frac{1}{\eps}\int_{0}^{1} \frac{d}{dt} \det D^{2} (u + t\eps \varphi_{\eps}) dt =
a^{ij}_{\eps}\partial_{ij}\varphi_{\eps}
\end{equation*}
where $$(a^{ij}_{\eps})_{ij}=\int_0^1 Cof \, \, ( D^{2}(u + t\eps \varphi_{\eps}))dt.$$ Because $u$ is strictly 
convex in 2D and $u_{\eps} \rightarrow u$
uniformly on $\overline B _1$, $D^{2}u_{\eps}\rightarrow D^{2} u$ uniformly on compact sets of $B_1$. Thus, as 
$\eps \to 0$, $a^{ij}_{\eps} \to U^{ij}$ uniformly on compact sets of $B_1$ and by Remark \ref{seq}, we find $\varphi_\eps \to \varphi$ uniformly in $\overline B_1$.
By the minimality of $u$, we find
$$0\le \lim_{\eps \to 0^{+}}\frac {1}{\eps}(L(u_\eps)-L(u))=\int_{\p B_1} \varphi \, d \sigma -\int_{B_1} \varphi \, dA.$$
By replacing $\varphi$ with $-\varphi$ we obtain the opposite inequality.
\end{proof}

\section{Proof of Proposition \ref{singular}}\label{qs}

In this section, we prove Proposition \ref{singular} where H1', H2' and H4' are satisfied. 
Given a convex function $u\in C^{\infty}(B_1)$ (not necessarily a minimizer of $L$) 
with $\rho \le \det D^2 u \le \rho^{-1}$, we let $v$ be the solution to the following Dirichlet problem
\begin{equation}\label{v}
U^{ij}v_{ij}=-A \quad \mbox{in $B_1$}, \quad \quad v=0 \quad \mbox{on $\p B_1$.}  
\end{equation}
Notice that $\Psi:=C(1-|x|^2)$ is an upper barrier for $v$ if $C$ is large enough, since
$$U^{ij} \Psi_{ij} \le -C\, tr \, U \le -C (\det D^2 U)^{1/n} = -C (\det D^{2}u)^{\frac{n-1}{n}}\leq
-C\rho^{\frac{n-1}{n}} \le -A,$$
hence
\begin{equation}\label{vbd}
0\le v(x) \le C(1-|x|^2) \, \, \sim \, dist(x, \p B_1).
\end{equation}
As in Lemma \ref{linear}, the function $v$ is the uniform 
limit of the corresponding $v_\delta$ that solve the Dirichlet problem in $B_{1-\delta}.$ Indeed, since $v_\delta$ also 
satisfies \eqref{vbd}, we see that $$|v_{\delta_1}-v_{\delta_2}|_{L^\infty} \le C \max\{\delta_1, \delta_2\}.$$

Let $\varphi$ be the solution of the homogenous problem
$$  U^{ij}\varphi_{ij}=0 \quad \mbox{in $B_1$}, \quad \quad \varphi=l^+ \quad \mbox{on $\p B_1$,}  $$
where $l^+=\max\{0, l\}$ for some linear function $l=b+\nu \cdot x$ of slope $|\nu|=1$. Denote 
by $\mathcal S:=\overline B_1 \cap \{l=0\}$ the segment of intersection of the crease of $l$ with $\overline B_1$. Then 

\begin{lem}\label{int_phi}
$$\int_{B_1}\varphi \, dA=\int_{B_1} l^+ \, dA +\int_{\mathcal S} u_{\tau \tau} v \, d \mathcal H^1,$$
where $\tau$ is the unit vector in the direction of $\mathcal S$, hence $\tau \perp \nu$.
\end{lem}

\begin{proof}
It suffices to show the equality in the case when $u\in C^{\infty}(\overline B_1)$. The general case follows by writing the identity in $B_{1-\delta}$ with $v_\delta$ (which increases as $\delta$ decreases),  and then letting $\delta \to 0$.

Let $\tilde l_\eps$ be a smooth approximation of $l^+$ with
$$D^2 l_\eps \rightharpoonup \nu \otimes \nu \, \, d H^1\lfloor \, {\mathcal S} \quad \mbox{as $\eps \to 0$,} $$
and let $\varphi_\eps$ solve the corresponding Dirichlet problem with boundary $\tilde l_\eps$.
Then, we integrate by parts and use $\partial_i U^{ij}=0,$
\begin{align*}
\int_{B_1}(\varphi_\eps-\tilde l_\eps) \, dA&=-\int_{B_1} (\varphi_\eps -\tilde l_\eps) U^{ij}v_{ij} \, dx\\
&=\int_{B_1} \partial_i(\varphi_\eps-\tilde l_\eps) U^{ij} v_j \, dx \\
&=-\int_{B_1} \partial_{ij}(\varphi_\eps-\tilde l_\eps) U^{ij} v \, dx \\
&=\int_{B_1} U^{ij} \partial_{ij}\tilde l_\eps \, \,  v \, dx.
\end{align*}
We let $\eps \to 0$ and obtain
$$\int_{B_1}(\varphi - l^+) \, dA=\int_{\mathcal S} U^{\nu \nu} v \, d \mathcal H^1,$$
which is the desired conclusion since $U^{\nu \nu}=u_{\tau \tau}$.
\end{proof}

From Lemma \ref{int_phi} and Proposition \ref{EL}, we obtain

\begin{cor}\label{u_tau_tau}
If $u$ is a minimizer of $L$ in the class $\mathcal A$ then
$$\int_{\mathcal S} u_{\tau \tau} v \, d \mathcal H^1=\int_{ \p B_1}l^+ \, d\sigma -\int_{B_1} l^+ \, dA.$$
The hypotheses on $\sigma$ and $A$ imply that if the segment $\mathcal S$ has length $2h$ with $h \le h_0$ small, universal then
$$ch^3 \le \int_{\mathcal S} u_{\tau \tau} v \, d \mathcal H^1 \le C h^3,$$
for some $c$, $C$ universal.
\end{cor}

\begin{lem}
 Let $X_{1}$ and $X_{2}$ be the endpoints of the segment $\mathcal S$ defined as above. Then
\begin{equation}
 \int_{\mathcal S} u_{\tau \tau} (1-\abs{x}^2)\, d \mathcal H^1 = 
4h \left(\frac{u(X_{1}) + u(X_{2})}{2}-\avint_{\mathcal S} u d \mathcal H^1\right),
\label{int-iden}
\end{equation}
where $2h$ denotes the length of $\mathcal S$.
\label{int-iden-lem}
\end{lem}
\begin{proof}
Again we may assume that $u \in C^2(\overline B_1)$ since the general case follows by approximating $B_1$ by $B_{1-\delta}$. Assume for simplicity that $\tau=e_1$. Then
$$\int_{\mathcal S} u_{\tau \tau} (1-\abs{x}^2) d \mathcal H^1 =\int_{-h}^{h} \partial_t^2 u(t,a) (h^2-t^2)\,dt$$
for some fixed $a$ and integrating by parts twice, we obtain (\ref{int-iden}).
\end{proof}

We remark that the right hand side in \eqref{int-iden} represents twice the area between the segment with end points $(X_1,u(X_1))$, $(X_2,u(X_2))$ and the graph of $u$ above $\mathcal S$.

\begin{defn}
We say that $u$ admits a tangent plane at a point $z \in \p B_1,$ if there exists a linear function $l_{z}$ such that
$$ x_{n+1}= l_{z}(x) $$ is a supporting
hyperplane for the graph of $u$ at $(z,u(z))$ but for any $\eps >0$,
$$ x_{n+1}=l_{z}(x)- \eps z \cdot (x-z) $$
is not a supporing hyperplane. We call $l_{z}$ a tangent plane for $u$ at $z$. 
\end{defn}

 \begin{rem}\label{dense}
 Notice that if $\det D^2 u \le C$ then the set of points where $u$ admits a tangent plane is dense in $\p B_1$. Indeed, using standard barriers it is not difficult to check that any point on $\p B_1$ where the boundary data $ u|_{\p B_1}$ admits a quadratic polynomial from below satisfies the definition above. In the definition above we assumed $u=\bar u$ on $\p B_1$ with $\bar u$ defined as in the Lemma \ref{l_semi}, therefore $ u|_{\p B_1}$ is lower semicontinuous. 
 \end{rem}
   
 Assume that $u$ admits a tangent plane at $z$, and denote by $$\tilde u=u-l_{z}.$$
 \begin{lem}\label{eta}
 There exists $\eta>0$ small universal such that the section
 $$\tilde S_{z}:=\{ x\in \overline B_1 | \quad \tilde u < \eta (x-z) \cdot (-z)\},$$
 satisfies $$\tilde S_{z} \subset B_1\setminus B_{1-\rho}, \quad |\tilde S_{z}| \ge c, $$
 for some small $c$ universal.
 \end{lem}  

\begin{proof}
We notice that \eqref{int-iden} is invariant under additions with linear functions. We apply it to 
$\tilde u$ with $X_1=z$, $X_2=x$ and use $\tilde u \ge 0$, $\tilde u(z)=0$ together with \eqref{vbd} and Corollary \ref{u_tau_tau} and obtain
$$\tilde u (x)\ge c|x-z|^2 \quad  x \in \p B_1 \cap B_{h_0}(z) .$$  
From the uniform strict convexity of $\tilde u$, which was obtained in Proposition \ref{strict-prop}, 
we find that the inequality above holds for all $x \in \p B_1$ for possibly a different value of $c$. Thus, by 
choosing $\eta$ sufficiently small, we obtain
$$\tilde S_{z} \subset B_1, \quad \quad \tilde S_{z} \cap B_{1-\rho}=\emptyset,$$ 
where the second statement follows also from Proposition \ref{strict-prop}.

Next we show that $|\tilde S_{z}|$ cannot be arbitrarily small. 
Otherwise, by the uniform strict convexity of $\tilde u$, we obtain that $\tilde S_{z} \subset B_{\eps^4}(z)$ for some small $\eps>0$.
 Assume for simplicity of notation that $z=-e_2$. Then the function
$$w:=\eta(x_2+1)+\frac \eps 2 x_1^2 + \frac {1}{2 \rho \eps} (x_2+1)^2 - 2\eps (x_2+1) ,$$
is a lower barrier for $\tilde u$ in $B_1 \cap B_{\eps^4}(z)$. Indeed, notice that if $\eps$ is sufficiently small then
$$w\le \eta(x_2+1) \le \tilde u \quad \mbox{on $\p (B_1 \cap B_{\eps^4} (z))$}, \quad \quad \det D^2 w =\rho^{-1}\ge \det D^2 \tilde u.$$ 
In conclusion, $\tilde u \ge w \ge (\eta/2)(x_2+1)$ and we contradict that $0$ is a tangent plane for $\tilde u$ at $z$.\end{proof}

\begin{lem}\label{sep_below}
Let $u$ be the normalized minimizer of $L$. Then $\|u\|_{C^{0,1}(\overline B_1)} \le C,$
and $u$ admits tangent planes at all points of $\p B_1$. Also, $u$ separates at least quadratically from its tangent planes i.e
$$u(x) \ge l_{z}(x)+c|x-z|^2, \quad \forall x,z \in \p B_1.$$
\end{lem}

\begin{proof}
Let $z$ be a point on $\p B_1$ where $u$ admits a tangent plane $l_{z}$. From the previous lemma we 
know that $u$ satisfies the quadratic separation inequality at $z$ and also that $\tilde u=u-l_{z}$ is bounded from 
above and below in $\tilde S_{z}$, i.e.,
$$|u-l_{z}| \le C \quad \mbox{in $\tilde S_{z}$}.$$
We obtain
$$
\int_{\tilde S_{z}} |l_{z}| \, dx -C  \le \int_{\tilde S_{z}} u \, d x \le \int_{B_1} u \, d x  
 \le C \int_{\p B_1} u \, d \sigma  \le C,
$$
and since $\tilde S_{z} \subset B_1$ has measure bounded from below we find $$l_{z}(z), |\nabla l_{z}| \le C.$$
By Remark \ref{dense}, this holds for a.e. $z \in \p B_1$ and, by approximation, we find that any point in $\p B_1$ admits a tangent plane that satisfies the bounds above. This also shows that $u$ is Lipschitz and the lemma is proved.  
\end{proof}

\begin{lem}\label{v_up_bd}
The function $v$ satisfies the lower bound
$$v(x) \ge c \, \, dist(x, \p B_1),$$
for some small $c$ universal.  
\end{lem}

\begin{proof}
Let $z \in \p B_1$ and let $l$ be a linear functional with
$$l(x)=l_z(x)-b \, z \cdot (x-z), \quad \mbox{for some $ 0 \le b \le \eta$}.$$
where  $l_z$ denotes a tangent plane at $z$.
We consider all sections $$S=\{x \in \overline B_1| \quad u<l\} $$ which satisfy
$$ \inf_{S}(u-l) \le -c_0,$$
for some appropriate $c_0$ small, universal. We denote the collection of such sections $\mathcal M_z$. From Lemma \ref{eta}, we
see that $\mathcal M_z \ne \emptyset$ since $\tilde S_z$ (or $b=\eta$) satisfies the property above. Notice also that $S \subset \tilde S_z \subset B_1$ and  $z \in \p S$. For any section $S\in \mathcal M_z$ we consider its center of mass $z^S$, and from the property above we see that $z^S \in B_{1-c}$ for some small $c>0$ universal. 

First, we show that the lower bound for $v$ holds on the segment $[z,z^S]$.  
Indeed, since 
$$U^{ij}[c(l-u)]_{ij} =-2c \det D^2 u \ge - 2c\rho^{-1} \ge -A =U^{ij}v_{ij},$$
and $c(l-u) \le 0 =v$ on $\p B_1$ we conclude that 
\begin{equation}\label{vbdbl}
c(l-u)^+ \le v \quad \mbox{in $B_1$}.
\end{equation} 
Now, we use the convexity of $u$ and the fact that the property of $S$ implies $(u-l)(z^S) <-c$, and conclude that
$$v(x) \ge c(l-u)(x) \ge c |x-z| \ge c \, \, dist(x, \p B_1) \quad \quad \forall x\in [z,z^S] .$$

Now, it remains to prove that the collection of segments $[z,z^S]$, $z \in \p B_1$, $S \in \mathcal M_z$ cover a fixed neighborhood of $\p B_1$. To this aim we show that the multivalued map
$$z \in \p B_1 \longmapsto F(z):=\{z^S| \quad S \in \mathcal M_z\} $$
has the following properties 

1) the map $F$ is {\it closed} in the sense that 
 $$z_n \to z_* \quad \mbox{and}  \quad z_n^{S_n} \to y_* \quad \Rightarrow y^*\in F(z_*)$$ 

2) $F(z)$ is a connected set for any $z$.

The first property follows easily from the following facts: $z^S$ varies continuously with the linear map $l$ that defines $S=\{u<l\}$; and if $l_{z_n} \to l_*$ then $l_* \le l_{z_*}$ for some tangent plane $l_{z_*}$. 

To prove the second property we notice that if we increase continuously the value of the parameter $b$ (which defines $l$) up to $\eta$ then all the corresponding sections belong also to $\mathcal M_z$. This means that in $F(z)$ we can connect continuously $z^S$ with $z^{\tilde S_z}$ for some section $\tilde S_z$. On the other hand the set of all possible $z^{\tilde S_z}$ is connected since the set $l_z$ of all tangent planes at $z$ is connected in the space of linear functions.

Since $F(z) \subset B_{1-c}$, it follows that for all $\delta<c$ the intersection map
$$z \longmapsto G_\delta (z)=\{ [z,y] \cap \p B_{1-\delta}| \quad y \in F(z)\}$$ 
has also the properties 1 and 2 above.  
Now it is easy to check that the image of $G_\delta$ covers the whole $\p B_{1-\delta}$, hence the collection of segments $[z,z^S]$ covers $B_1\setminus B_{1-c}$ and the lemma is proved.

\end{proof}

Now, we are ready to prove the first part of Proposition \ref{singular}.

\begin{proof}[Proof of Proposition \ref{singular} (i)]

In Lemma \ref{eta}, we obtained the quadratic separation from below for $\tilde u=u-l_{z}$. Next we show that $\tilde u$ separates at most quadratically on $\p B_1$ in a neighborhood of $z$.

 Assume for simplicity of notation that $z=-e_2$. We apply \eqref{int-iden} to $\tilde u$ with $X_1=(-h,a)$, $X_2=(h,a)$, then use Corollary \ref{u_tau_tau} and Lemma \ref{v_up_bd} and obtain
$$\frac{\tilde u(X_{1}) +\tilde u(X_{2})}{2}-\avint_{\mathcal S} \tilde u \le C h^2.
 $$
 On the other hand, for small $h$, the segment $[z,z^{\tilde S_z}]$ intersects $[X_1,X_2]$ at a point
 $y=(t,a)$ with $|t| \le C h^2 \le h/2$. Moreover, since $y \in \tilde S_z$ we have
 $\tilde u(y) \le \eta (a+1) \le C h^2$. On the segment $[X_1,X_2]$ , $\tilde u$ satisfies the conditions of Lemma \ref{1D_lem} which we prove below, hence
 $$\tilde u(X_1), \tilde u(X_2) \le C h^2.$$
 In conclusion, $u$ separates quadratically on $\p B_1$ from its tangent planes and therefore satisfies the 
hypotheses of the Localization Theorem in \cite{S2}, \cite{LS}. From Theorem 2.4 and Proposition 2.6 in \cite{LS}, we conclude that 
 \begin{equation}\label{beta}
 \|u\|_ {C^{1,\beta}(\overline B_1)},\, \, \|v\|_{C^\beta(\overline B_1)} ,\, \, \|v_\nu\|_{C^\beta(\p B_1)}\le C, 
 \end{equation}
 for some $\beta<1$, $C$ universal.
 \end{proof}

\begin{lem}\label{1D_lem}
 Let $f:[-h,h]\to \R^+$ be a nonnegative convex function such that
\begin{equation*}
 \frac{f (-h) + f(h)}{2}-\frac{1}{2h}\int_{-h}^{h} f(x) dx \leq Mh^2, \quad \quad f(t) \le Mh^2,
\end{equation*}
for some $t\in [-h/2, h/2]$. Then
\begin{equation*}
 f(\pm h)\leq Ch^2,
\end{equation*}
for some $C$ depending on $M$.
\end{lem}

\begin{proof}
The inequality above states that the area between the line segment with end points $(-h,f(-h))$, $(h,f(h))$ and the graph of $f$ is bounded by $2M h^3$. By convexity, this area is greater than the area of the triangle with vertices $(-h,f(-h))$, $(t,f(t))$, $(h,f(h))$. Now the inequality of the heights $f(\pm h)$ follows from elementary euclidean geometry.   
\end{proof}

Finally, we are ready to prove the second part of Proposition \ref{singular}.

\begin{proof}[Proof of Proposition \ref{singular} (ii)]
Let 
$\varphi$ be such that
$$U^{ij} \varphi_{ij} =0 \h~\text{in} ~\Omega, \quad \quad \varphi \in C^{1,1}(\p B_1) \cap C^0(\overline B_1).$$
Since $u$ satisfies the quadratic separation assumption and $f$ is smooth up to the boundary, we obtain from Theorem 2.5 and Proposition 2.6 in \cite{LS}
$$\|v\|_{C^{1, \beta}(\overline B_1)}, \, \,   \|\varphi\|_{C^{1, \beta}(\overline B_1)} \le K , \quad \quad \mbox{and} \quad |U^{ij}| \leq K |\log \delta|^2 \quad \mbox{on $B_{1-\delta,}$}  $$
for some constant $K$ depending on $\rho$, $\|f\|_{C^\beta(\overline B_1)},$ and $\|\varphi\|_{C^{1,1}(\p B_1)}$. 

We will use the following identity in 2D:
\begin{equation*}
 U^{ij} v_{j}\nu_{i} = U^{\tau\nu} v_\tau + U^{\nu\nu} v_{\nu}.
\end{equation*}
Integrating by parts twice, we obtain as in \eqref{intparts} 
\begin{align*}
\int_{B_{1-\delta}} \varphi \, dA&=-\int_{B_{1-\delta}}  \varphi \, U^{ij}v_{ij} \, dx\\
&=\int_{\p B_{1-\delta}} \varphi_i  U^{ij}v \nu_j -\int_{\p B_{1-\delta}} \varphi U^{ij}v_j \nu_i\\
&=-\int_{\p B_{1-\delta}} \varphi U^{\nu \nu} v_\nu + o(\delta)
\end{align*}
where in the last equality we used the estimates $$|v| \le C \delta,  \quad |v_\tau| \le K \delta^\beta, \quad |\varphi|, |\nabla \varphi| \le K, \quad U^{ij} \le K |\log \delta|^2 \quad \mbox{on $\p B_{1-\delta}$.}$$ 

Since on $\p B_r$ $$U^{\nu \nu}=u_{\tau \tau}= r^{-2} u_{\theta \theta} + r^{-1} u_\nu,$$
$u \in C^{1, \beta}(\overline B_1)$ and $u(r e^{i \theta})$ converges uniformly as $r \to 1$, and $u_{\theta \theta}$ is 
uniformly bounded from below, we obtain
$$U^{\nu \nu} \, d \mathcal H^1 \lfloor _{ \p B_r} \rightharpoonup (u_{\theta \theta} + u_\nu) \,  d \mathcal H^1 \lfloor_{\p B_1}  \quad \mbox {as $r \to 1$}.$$
We let $\delta \to 0$ in the equality above and find
$$ \int_{B_1} \varphi \, dA=-\int_{\p B_1} \varphi \, (u_{\theta \theta}+u_\nu) v_\nu \, d \mathcal H^1.$$
Now the Euler-Lagrange equation, Lemma \ref{EL}, gives
$$ (u_{\theta \theta}+u_\nu) v_\nu=-\sigma \quad \mbox{on $\p B_1$.}  $$
We use that  $\|v_\nu\|_{C^\beta(\p B_1)} \le C$ and, from Lemma \ref{v_up_bd}, $v_\nu \le -c$ on $\p B_1$ and obtain $$\|u\|_{C^{2,\gamma}(\p B_1)} \le C \|\sigma\|_{C^\gamma(\p B_1)}.$$ 
\end{proof}

\section{The general case for $A$} \label{5}

In this section, we remove the assumptions that $A$ is bounded from 
below by $\rho$ in $B_1$ and also we assume that $A$ is bounded from above only in a neighborhood of the boundary. Precisely, 
we assume that  $A \ge 0$ in $B_1$ and $A \le \rho^{-1}$ in $B_{1}\backslash \overline B_{1-\rho}$. We may also 
assume $A$ is smooth in $B_1$ since the general case follows by approximation. Notice that $\int_{B_1} A \, dx$ is bounded 
from above and below since it equals  $\int_{\p B_1} d \sigma$.

Let $v$ be the solution of the Dirichlet problem
\begin{equation}\label{v2}
U^{ij}v_{ij}=-A, \quad \quad v=0 \quad \mbox{on $\p B_1$.}
\end{equation}

In Section \ref{qs}, we used that $A$ is bounded from above when we obtained $v \le C (1-|x|^2)$, and we used that $A$ is 
bounded from below in Lemma \ref{v_up_bd} (see \eqref{vbdbl}). We need to show that these bounds for $v$ also hold in a 
neighborhood of $\p B_1$ under the weaker hypotheses above. First, we show 
\begin{lem}\label{genA}
$$v \le C  \quad \mbox{on $\p B_{1-\rho/2}$,} \quad \quad 
v \ge c(\delta) \quad \mbox{on $B_{1-\delta}$,} $$ 
with $C$ universal, and $c(\delta)>0$ depending also on $\delta$. 
\end{lem}

\begin{proof} As before, we may assume that $u \in C^{\infty}({\overline B_1})$ since the general case follows by approximating $B_1$ by $B_{1-\eps}$.

We multiply the equation in \eqref{v2} by $(1-|x|^2)$, integrate by parts twice and obtain
$$\int _{B_1} 2 v \, \, tr \, U \, dx=\int_{B_1} A(x) (1-|x|^2) \, dx \le C,$$
and since $tr \, U \ge c$ we obtain $$\int_{B_1} v \, dx \le C.$$
We know 

1) $v \ge 0$ solves a linearized Monge-Amp\`ere equation with bounded right hand side in $B_1 \setminus B_{1-\rho}$ ,

2) $u$ has a uniform modulus of convexity on compact sets of $B_1$.

\noindent
Now we use the Harnack inequality of Caffarelli-Gutierrez \cite{CG} and conclude that
$$ \sup_{\mathcal V} v \le C (\inf_{\mathcal V} v +1),\quad \quad \mathcal V:=B_{1-\rho /4} \setminus \overline B_{1-3 \rho/4},$$ and the integral inequality above gives $\sup_{\mathcal V} v \le C$.

Next, we prove the lower bound. We multiply the equation 
in \eqref{v2} by $\varphi \in C_0^\infty(B_1)$ 
with $$\varphi=0 \quad \mbox{if $|x| \ge 1-\delta/2$}, \quad \varphi=1 \quad \mbox{in $B_{1- \delta}$}, 
\quad \|D^2 \varphi \| \le C/ \delta^2,  $$ integrate by parts twice and obtain
$$ C(\delta) \int_{\mathcal U} v \, \, tr \, U  \ge - \int_{B_1} v \, U^{ij} \varphi_{ij} \, = \int_{B_1} A \varphi  \ge c, 
\quad \quad \mathcal U:=B_{1-\delta/2} \setminus \overline B_{1-\delta},$$
where the last inequality holds provided that $\delta$ is sufficiently small. Since $u$ is normalized we 
obtain (see Proposition \ref{strict-prop}) , $|\nabla u| \le C(\delta)$ in $\mathcal U$ thus
$$\int_{ {\mathcal U}} tr \, U =\int_{ {\mathcal U}} \triangle  u  = \int_{\p  {\mathcal U}}u_\nu  \le C(\delta).$$
The last two inequalities imply $\sup_{{\mathcal U}}v \ge c(\delta)$, hence there exists $x_0 \in {\mathcal U}$ such 
that $v(x_0) \ge c(\delta).$ We use 1), 2) above and Harnack inequality and find 
$v \ge c(\delta)$ in $B_{\bar \delta}(x_0)$ for some small $\bar \delta$ depending on $\rho$ and $\delta$. Since $v$ is a 
supersolution, i.e $U^{ij}v_{ij} \le 0$, we can apply the weak Harnack inequality of Caffarelli-Gutierrez, Theorem 4 in \cite{CG}. From property 2) above, we see that we can extend the lower bound of $v$ from 
$B_{\bar \delta}(x_0)$ all the way to $\mathcal U$, and by the maximum principle this bound holds also in $B_{1-\delta/2}$. 

\end{proof}

The upper bound in Lemma \ref{genA} gives as in \eqref{vbd} the upper bound for $v$ in a neighborhood of $\p B_1$, i.e
$$v(x) \le C (1-|x|^2) \quad \mbox{on} \quad B_1 \setminus B_{1-\rho/2}.$$
This implies as in Section \ref{qs} that Lemma \ref{sep_below} holds i.e., $u$ separates at least quadratically 
from its tangent planes on $\p B_1$. It remains to show that also Lemma \ref{v_up_bd} holds. Since $A$ is not 
strictly positive, $c(l-u)$ is no longer a subsolution for the equation \eqref{v2} and we cannot 
bound $v$ below as we did in \eqref{vbdbl}. In the next lemma, we construct another barrier which alows us to bound $v$ from below on the segment $[z,z^S]$.

\begin{lem}\label{hopf}
Let $\tilde u: B_1 \to \R$ be a convex function with $\tilde u \in C(\overline B_1) \cap C^2(B_1)$, and 
$$\rho \le \det D^2 \tilde u \le \rho^{-1}.$$ Assume that the section $S:=\{ \tilde u<0\}$ is included in $B_1$ and is 
tangent to $\p B_1$ at a point $z \in \p B_1$, and also that
$$\inf_S \tilde u \le - \mu,$$ for some $\mu>0$. If 
$$\tilde U ^{ij} v_{ij} \le 0 \quad \mbox{in $B_1$}, \quad v \ge 0 \quad \mbox{on $\p B_1$,}$$
then $$v(x) \ge c(\mu, \rho) |x-z| \, \, \inf_{S'} v\quad \forall x\in [z,z^S],\quad \quad S':=\{\tilde u \le \frac 12 \inf_S \tilde u \},$$
where $z^S$ denotes the center of mass of $S$, and $c(\mu,\rho)$ is a positive constant depending on $\mu$ and $\rho$.
\end{lem} 

The functions $\tilde u= u-l$ and $v$ in the proof of Lemma \ref{v_up_bd}  satisfy the lemma above, if $\eta$ in Lemma
\ref{eta} is small, universal. Using also the lower bound on $v$ from Lemma \ref{genA}, we find
$$v \ge c|x-z| \quad \mbox{on $[z,z^S]$},$$
for some $c$ universal, and the rest of the proof of Lemma \ref{v_up_bd} follows as before. This shows that Proposition \ref{singular} holds also with our assumptions on the measure $A$. 

\begin{proof}[Proof of Lemma \ref{hopf}]
We construct a lower barrier for $v$ of the type
$$ w:=e^{k \bar w}-1, \quad \bar w:=-\tilde u + \frac \eps 2 (|x|^2-1),$$
for appropriate constants $k$ large and $\eps \ll \mu$ small. Notice 
that $w \le 0$ on $\p B_1$ since $\bar w \le 0$ on $\p B_1$.  Also $$\bar w \ge c \,  |x-z|  \quad \mbox{on $[z,z^S]$},$$
since, by convexity, $-\tilde u \ge c|x-z|$ on $[z,z^S]$ for some $c$ depending on $\mu$ and $\rho$.
It suffices to check that
$$\tilde U^{ij} w_{ij} \ge 0 \quad \mbox{on $B_1 \setminus S'$,}$$ since then we obtain $v \ge (\inf_{S'}v) \, c  \, w$ in $B_1 \setminus S'$ which easily implies the conclusion.
In $B_1 \setminus S'$ we have $|\nabla \bar w| \ge c(\mu)>0$ provided that $\eps$ is sufficiently small, thus $$\tilde U ^{ij} \bar w_i \bar w_j =(\det D^2 \tilde u)  (\nabla \bar w)^T \, (D^2 \tilde u)^{-1}  \, \nabla \bar w \ge  c \Lambda^{-1},$$ where $\Lambda$ is the largest eigenvalue of $D^2 \tilde u$. Then, we use that $tr \, \tilde U \ge c \lambda^{-1}\ge c \Lambda^\frac {1} {n-1}$ where $\lambda$ is the smallest eigenvalue of $D^2 \tilde u$, and obtain
\begin{align*}
\tilde U^{ij} w_{ij}&=ke^{k \bar w}\left (\tilde U^{ij}\bar w_{ij} +k \tilde U^{ij} \bar w_i \bar w_j  \right ) \\
&\ge ke^{k \bar w}\left (-n + \eps \, tr \, \tilde U + k c \Lambda^{-1}  \right )\\
&\ge  ke^{k \bar w}\left (-n + c(\eps \, \Lambda^\frac{1}{n-1}  + k  \Lambda^{-1})  \right )\\
&\ge 0,
\end{align*}
if $k$ is chosen large depending on $\eps$, $\rho$, $\mu$ and $n$.
\end{proof}

\section{Singular minimizers in dimension $n \ge 3$.} \label{Sec-sing}

Let $$u(x):=|x'|^{2-\frac 2 n} h(x_n),$$ be the singular solution to $\det D^2u =1$ constructed by Pogorelov, with $h$ a smooth even function, defined in a neighborhood of $0$ and $h(0)=1$, satisfying an ODE $$\left((1-\frac2 n)hh''-(2-\frac 2 n)h'^2\right )h^{n-2}=c.$$  
We let
$$v(x):=|x'|^{2-\frac 2 n} q(x_n)$$
be obtained as the infinitesimal difference between $u$ and a rescaling of $u$, $$v(x',x_n):=\lim_{\eps \to 0} \frac 1 \eps [u(x',x_n)-(1+\eps)^{-\gamma} u(x',(1+\eps)x_n)],$$
for some small $\gamma < 2/n$.
Notice that $$q(t)= \gamma h(t)-h'(t)t $$ and $q>0$ in a small interval $(-a,a)$ and $q$ vanishes at its end points. Also, $$U^{ij} v_{ij}= n \gamma -2 <0 \quad \mbox{in} \quad 
\Omega:= \R^{n-1} \times  [-a,a],$$
$$v=0, \quad U^{\nu \nu} v_{\nu}=U^{nn}v_n=-\sigma_0 \quad \mbox{on} \quad \p \Omega,$$
for some constant $\sigma_0>0$. The last equality follows since $U^{nn}$ is homogenous of degree $-(n-1)(2/n)$ in $|x'|$ and $v_n$ is homogenous of degree $2-2/n$ in $|x'|$.

Notice that if $u$, $v$ are solutions of the system \eqref{EL-eq-intro} in the infinite cylinder $\Omega$ for uniform measures $A$ and $\sigma$. In order to obtain a solution in a finite domain $ \Omega_0$ we modify $v$ outside a neighborhood of the line $|x'|=0$ by subtracting a smooth convex function $\psi$ which vanishes in $B_1$ and increases rapidly outside $B_1$. Precisely we let
$$\tilde v:=v-\psi, \quad  \Omega_0:= \{ \tilde v >0\}$$ and then we notice that $u$, $\tilde v$,  solve the system \eqref{EL-eq-intro} in the smooth bounded domain $ \Omega_0$ for smooth measures $A$ and $\sigma$. 
 
 Since $$|U^{ij}| \le C r^{\frac 2n-2} , \quad \mbox{if $|x'| \ge r$},$$ we integrate by parts in the domain $ \Omega_0 \setminus \{|x'| \le \eps\}$ and then let $\eps \to 0$ and find
 $$\int_{ \Omega_0} \varphi \, dA=-\int_{\Omega_0} U^{ij}\varphi_{ij} v +\int_{\p \Omega_0} \varphi \, d \sigma, \quad \quad \forall \varphi \in C^2(\overline  \Omega_0),$$
or $$L (\varphi)=\int_{\Omega_0} U^{ij}\varphi_{ij} \, v.$$

This implies that $L$ is stable, i.e $L(\varphi)>0$ for any convex $\varphi$ which is not linear. Also, if $w \in C^2(\overline \Omega_0)$ satisfies $\det D^2 w=1$, then $U^{ij}(w-u)_{ij} \ge 0$,  and we obtain
$$L(w)-L(u)=\int_{\Omega_0}U^{ij}(w-u)_{ij} \,  v \ge 0,$$
i.e $u$ is a minimizer of $L$. 

We remark that the domain $\Omega_0$ has flat boundary in a neighborhood of the line $\{ |x'|=0\}$ and therefore is not uniformly convex. However this is not essential in our example. One can construct for example a function $\bar v$ in a uniformly convex domain by modifying $v$ as
$$\bar v:= |x'|^{2-\frac 2n} q(x_n(1+\delta |x'|^2)),$$
for some small $\delta>0$.

\section{Proof of Theorem \ref{min_func}}\label{last}

We assume for simplicity that $\Omega=B_1$. The existence of a minimizer $u$ for the convex functional $E$ follows as in Section \ref{Sec2}. 
First, we show that
\begin{equation}\label{det_bd}
t_1 \le \det D^2u \le t_0
\end{equation}
for some $t_1$ depending on $F$ and $\rho$. The upper bound follows easily. If $\det D^2u >t_0$ in a set of positive measure then the function $w$ defined as
$$\det D^2 w=\min \{ \, \, t_0, \, \, \det D^2 u \} , \quad \quad \mbox{$w=u$ on $\p B_1$},$$
satisfies $E(w)<E(u)$ since $F(\det D^2 w)=F(\det D^2u)$ and $L(w)<L(u)$.

In order to obtain the lower bound in \eqref{det_bd} we need the following lemma. 
\begin{lem}\label{alexandrov}
Let $w$ be a convex functions in $ B_1$ with $$(\det D^2 w)^{\frac 1n}=g \in L^n(B_1).$$
Let $w+\varphi$ be another convex function in $B_1$ with the same boundary values as $w$ such that 
$$(\det D^2(w+\varphi))^{\frac 1n} =g-h, \quad \mbox{for some $h\ge 0$}.$$
 Then
$$\int_{B_1}\varphi \, g^{n-1} \le C(n) \int_{B_1}h \, g^{n-1} .$$
\end{lem}

\begin{proof}
By approximation, we may assume that $w$, $\varphi$ are smooth in $\overline B_1$. Using the concavity of the map $M \mapsto (\det M)^{\frac 1n}$ in the space of symmetric matrices $M \ge 0$, we obtain
$$(\det D^2(w+\varphi)) ^{\frac 1 n} \le (\det D^2 w)^{\frac 1n} +\frac 1n (\det D^2 w)^{\frac 1n -1} W^{ij} \varphi_{ij},$$
hence
$$-n \, h \, g^{n-1} \le W^{ij} \varphi_{ij}.$$
We multiply both sides with $\Phi:=\frac 12 (1-|x|^2)$ and integrate. Since both $\varphi$ and $\Phi$ vanish on $\p B_1$ we integrate by parts twice and obtain
$$- C(n) \int_{B_1} h \, g^{n-1} \le \int_{B_1} W^{ij} \Phi_{ij} \, \, \varphi =-\int_{B_1} (tr \, \, W) \varphi.$$
Using $$tr \, W \ge c(n) (\det W)^\frac 1n=c(n) (\det D^2w)^\frac{n-1}{n} =c(n) g^{n-1}$$
we obtain the desired conclusion.
\end{proof}

Now we prove the lower bound in \eqref{det_bd}. Define $w$ such that $w=u$ on $\p B_1$ and $$\det D^2 w=\max \{\, \, t_1, \, \det D^2u  \},$$
for some small $t_1$. Since $G(t)=F(t^n)$ is convex and $\det D^2 w \geq t_{1},$ we have    
 $$G((\det D^2 w)^{1/n}) \le G((\det D^2 u)^{1/n}) + G'(t_1^{1/n}) ((\det D^2w)^{1/n}-(\det D^2u)^{1/n}) .$$ We denote $$u-w=\varphi, \quad (\det D^2 w)^{1/n}=g, \quad (\det D^2 u)^{1/n}=g-h,$$and we rewrite the inequality above as
$$F(\det D^2 w) \le F(\det D^2 u) + G'(t_1^{1/n}) \, h.$$ From Lemma \ref{alexandrov}, we obtain
$$\int_{B_1} h \, g^{n-1} \ge c(n) \int_{B_1} \varphi \,g^{n-1}$$ and since $h$ is supported on the set where the 
value of $g=t_1^{1/n}$ is minimal, we find that 
$$\int_{B_1} h \ge c(n) \int_{B_1} \varphi .$$ This gives
$$\int_{B_1} F(\det D^2 w)-F(\det D^2 u) \le c(n) G'(t_1^{1/n}) \int_{B_1}\varphi,$$ thus, using the minimality of $u$ and $G'(0^+)=-\infty$,
$$0 \le E(w)-E(u) \le \int_{B_1} \varphi dA +  c(n) G'(t_1^{1/n}) \int_{B_1}\varphi \le 0, $$
if $t_1$ is small enough. In conclusion, $\varphi=0$ and $u=w$ and \eqref{det_bd} is proved.

\

We denote
$$ \det D^2 u=f, \quad t_1 \le f \le t_0.$$
Any minimizer for $L$ in the class of functions whose determinant equals $f$ is a minimizer for $E$ as well. In order to 
apply Theorem \ref{C2-thm} we need $f$ to be Holder continuous. However, we can approximate $f$ by smooth functions $f_n$ and find smooth minimizers $u_n$ for 
approximate linear functionals $L_n$ with the constraint $\det D^2 u_n=f_n$. By Proposition \ref{singular} (see \eqref{beta}), 
$$\|u_n\|_{C^{1,\beta}(\overline B_1)}, \, \, \|v_n\|_{C^\beta(\overline B_1)} \, \le C,$$ hence we may 
assume (see Theorem \ref{cpt}) that, after passing to a subsequence, $u_n \to u$ and $v_n \to v$ uniformly for 
some function $v \in C^\beta(\overline B_1)$. We show that 
\begin{equation}\label{v_f}
v=-F'(f). 
\end{equation}
Then by the hypotheses on $F$ we obtain $\det D^2u=f \in C^\beta(\overline B_1)$ and from Theorem \ref{C2-thm} we easily obtain $$\|u\|_{C^{2,\alpha}(\overline B_1)}, \, \, \|v\|_{C^{2,\alpha}(\overline B_1)} \le C,$$ for some $C$ depending on $\rho$, $\alpha$, $\|\sigma\|_{C^\alpha(\overline B_1)}$, $\|A\|_{C^\alpha(\overline B_1)}$ and $F$. 

In order to prove \eqref{v_f} we need a uniform integral bound (in 2D) between solutions to the Monge-Amp\`ere equation and solutions of the corresponding linearized equation.

\begin{lem}\label{lin}
Assume $n=2$ and let $w$ be a smooth convex function in $\overline B_1$ with
$$\lambda \le \det D^2 w:=g \le \Lambda,$$ for some positive constants $\lambda$, $\Lambda$. 
Let $w+\eps \varphi$ be a convex function with  
$$\det D^2 (w+ \eps \varphi)=g+\eps h, \quad \varphi=0 \quad \mbox{on $\p B_1$} $$
for some smooth function $h$ with $\|h\|_{L^\infty}\le 1$. If $\eps \le \eps_0$ then
$$\int_{B_1}|h-W^{ij} \varphi_{ij}| \le C \eps.$$
for some $C$, $\eps_0$ depending only on $\lambda$, $\Lambda$.
\end{lem}

We postpone the proof of the lemma untill the end of the section.

Now let $h$ be a smooth function, $\|h\|_{L^\infty}\le 1$, and we solve the equations 
$$\det D^2(u_n+\eps \varphi_n)=f_n+\eps h, \quad \varphi_n=0 \quad \mbox{on $\p B_1$},$$
with $u_n$, $f_n$ as above. From \eqref{intparts} we see that 
$$L_n(\varphi_n)=\int_{B_1} (U_n^{ij} \p_{ij}\varphi_n) \, v_n,$$ 
hence, by the lemma above
$$|L_n(\varphi_n)-\int_{B_1}h \, v_n| \le C \eps$$
with $C$ universal. We let $n \to \infty$ and obtain 
$$|L(\varphi)-\int_{B_1} h\, v| \le C \eps.$$
with $\varphi$ the solution of
$$\det D^2(u+\eps \varphi)=f+\eps h, \quad \varphi=0 \quad \mbox{on $\p B_1$}.$$

The inequality $E(u+\eps \varphi) \ge E(u)$ implies
$$\int_{B_1}(F(f+\eps h)-F(f) +  \eps h \, v) \ge -C \eps^2,$$ hence, as $\eps \to 0$,
$$\int_{B_1} (F'(f)+v) \, h \ge 0 \quad \mbox{for any smooth $h$},$$
which gives \eqref{v_f}.
\qed

\begin{proof}[Proof of Lemma \ref{lin}]
Using the concavity of $(\det D^2 w)^{1/n}$ we obtain
$$(g+\eps h)^{1/n} \le g^{1/n} + \frac \eps n g^{1/n -1} W^{ij} \varphi_{ij} ,$$
thus, for $\eps\leq \eps_{0}$
\begin{equation}\label{one_side}
h- C \eps \le W^{ij} \varphi_{ij}.
\end{equation}
Since $n=2$ we have
$$\det D^2(w+\eps \varphi)=\det D^2 w + \eps W^{ij} \varphi_{ij}+ \eps^2 \det D^2 \varphi,$$ hence
$$h- W^{ij}\varphi_{ij}=\eps \det D^2 \varphi.$$
From the pointwise inequality \eqref{one_side}, we see that in order to prove the lemma it suffices to show that
$$\int_{B_1}\det D^2 \varphi \ge -C.$$

Integrating by parts and using $\varphi=0$ on $\p B_1$ we find
$$\int_{B_1}2\det D^2 \varphi=\int_{B_1}\Phi^{ij}\varphi_{ij}=\int_{\p B_1} \Phi^{ij} \varphi_i \nu_j=\int_{\p B_1} \Phi^{\nu \nu} \varphi_\nu=\int_{\p B_1}\varphi_\nu^2 \ge 0$$
where we used that $\Phi^{\nu \nu}=\varphi_{\tau \tau}=\varphi_\nu$.
\end{proof}

\end{document}